\documentclass[english]{cccconf}
\usepackage[comma,numbers,square,sort&compress]{natbib}
\usepackage{empheq,subeqnarray,enumerate,color}
\usepackage{mathrsfs,dsfont,amssymb,amsthm,amsmath,amsfonts,bbm,graphicx,cases}

\newtheorem{condition**}{A*}
\newtheorem{condition***}{C*}
\newtheorem{condition*}{C}

\newtheorem{example}{Example}[section]

\newtheorem{theorem}{Theorem}[section]
\newtheorem{lemma}{Lemma}[section]

\begin{document}

\title{Partial Information Differential Games for Mean-Field SDEs}
\author{Hua XIAO\aref{sdu},
        Shuaiqi ZHANG\aref{GDIU},
        %WANG Wu\aref{hit}
        }

% Note: the first argument in the \affiliation command is optional.
% It defines a label for the affiliation which can be used in the \aref
% command. If there is only one affiliation for all authors, then the
% optional argument in the \affiliation command should be suppressed,
% and the \aref command should aslo be removed after each author in
% \author command, in this case the affiliation will not be numbered.

% Çë×¢Ò⣺\affiliationÃüÁîµÄµÚÒ»¸ö²ÎÊýÊÇ¿ÉÑ¡µÄ£¬Ëü¶¨ÒåÁËÓÃÓÚ\arefÃüÁîµÄ±êÇ©¡£
% Èç¹ûËùÓÐ×÷ÕßÖ»ÓÐÒ»¸öµ¥Î»£¬Ç벻ҪʹÓÃ\affiliationÃüÁîµÄ¿ÉÑ¡²ÎÊý£¬Í¬Ê±ÔÚÉÏÃæ
% µÄ\authorÃüÁîÖеÄÿλ×÷ÕßÐÕÃûºóÃæÒ²²»ÄÜʹÓÃ\arefÃüÁʾÀýÈçÏÂ
% \author{First Author, Second Author, Third Author}
% \affiliation{Chinese Academy of Sciences, Beijing 100190, P.~R.~China\email{ccc@amss.ac.cn}}
% ´Ëʱµ¥Î»Ç°²»»áÓÐÊý×Ö±àºÅ£¬×÷ÕßÐÕÃûºóÃæҲûÓбàºÅ

\affiliation[sdu]{School of Mathematics and Statistics, Shandong University, Weihai 264209, China
        \email{xiao$\_$hua@sdu.edu.cn}}
\affiliation[GDIU]{School of Economics and Commerce, Guangdong University of Technology, Guangzhou 510520, China \email{shuaiqiz@hotmail.com}}

\maketitle

\begin{abstract}
This paper is concerned with non-zero sum differential games of mean-field stochastic differential equations with partial information and convex control domain. First, applying the classical convex variations, we obtain stochastic maximum principle for Nash equilibrium points. Subsequently, under additional assumptions,  verification theorem for Nash equilibrium points is also derived. Finally, as an application, a linear quadratic example is discussed. The unique Nash equilibrium point is represented in a feedback form of not only the optimal filtering but also expected value of the system state, throughout the solutions of the Riccati equations.
\end{abstract}


% Please remove or comment out the following line if the footnote is not necessary
\footnotetext{This work is partially supported by National Natural Science
Foundation (NNSF) of China under Grant 11371228, 11471192, 61573217 and 11501129.}

\section{Introduction}% Section 1.

In this paper, we study partial information stochastic differential game problems in which system states are governed by stochastic differential equations (SDEs) of mean-field type, in the sense that the coefficients of the SDEs depend not only on the system states, but also on their expected values.
Also, the SDEs of mean-field type are often used to describe the aggregate behavior of lots of mutually interacting particles at mesoscopic level and play an important role in physics, finance, economics, etc. For more information,  we refer the reader, for instance, to \cite{Chan1994,Kac1958} as well as the references therein. Recently, a new kind of backward SDEs (BSDEs) of mean-field type has been studied by Buckdahn et al. \cite{BDLP2009,BLP2009} which is called mean-field BSDEs. For classical control problems of SDEs without mean field, we refer the readers to \cite{Peng1990,YZ1999}, etc.

Mean-field games and mean-field control problems have received considerable attention in the probability and optimal control literature in recent years. Li \cite{Li2012} studied the stochastic maximum principle for mean-field SDEs with convex control domain and also got the verification theorem under additional conditions. Buckdahn, Djehiche and Li \cite{BDL2011} used the classical spike perturbation and derived a Peng-type general stochastic maximum principle. Yong \cite{Yong2013} investigated linear-quadratic (LQ) optimal control problems for mean-field SDEs and a feedback representation was obtained for the optimal control.  Lasry and Lions \cite{LL2007} presented three examples of
mean-field approach to modelling in Economics and Finance, derived nonlinear mean-field SDEs and established their links with various fields of Analysis. More recent developments and their applications of mean-fiels games of SDEs can be found in Bensoussan, Sung and Yam \cite{BSY2013}, Carmona, Delarue and Lachapelle \cite{CDL2013}, Gu¨¦ant \cite{Gu2009}, etc., and the references therein. Different from the above works, we consider two players non-zero sum differential games of mean-field SDEs with partial information and convex control domain. The distinguishing feature is the information available to the two players is the sub-filtration of full information. The problem we study may cover many control and game problems of mean-field SDEs with complete information as special cases. The present work will also enrich the relevant theory of stochastic filtering.

The rest of this paper is organized as follows. In Section 2, we specify the problem considered. Section 3 is devoted to deriving the stochastic maximum principle and verification theorem for Nash equilibrium points. Finally, in Section 4, we solve an LQ example to explain our application. By introducing the systems of some Riccati equations and forward-backward stochastic filtering equations of mean-filed type, we give the feedback representation for the unique Nash equilibrium point.

\section{Formulation of Problem}

Let $|x|$ denote the Euclidian norm of $x\in\mathbb{R}^{n}$ and
$\langle x, y\rangle$ be the inner product of $x,
y\in\mathbb{R}^{n}$. The transpose and Euclidian norm of a matrix
$M=(m^{ij})_{{1\leq i\leq n \atop 1\leq j\leq d}}=(m^{1}, \cdots,
m^{d})\in\mathbb{R}^{n\times d}$ are expressed as $M^{*}$ and
$|M|=\sqrt{\hbox{trace}(MM^{*})}$, respectively. Similarly, $\langle
M_1, M_2 \rangle$ = trace $(M_1M_2^{*})$ with $M_1,
M_2\in\mathbb{R}^{n\times d}$.
Let $T>0$ be a fixed constant and $C$ be a positive constant which can be different from line to line. Let $\left(\Omega, \mathscr{F}, (\mathscr{F}_{t})_{0\leq t\leq T}, \mathbb{P}\right)$ be a complete filtered probability space on which $\mathscr{F}_{t}$ denotes a natural filtration generated by a standard Brownian motion $(w_{1}, w_{2})$ with values in $\mathbb{R}^{d_{1}+d_{2}}$.

We only consider the case of two players and define the admissible control set $\mathscr{U}_{i}$ for Player $i\; (i=1, 2)$ by
\begin{multline}\label{Ad}
    \mathscr{U}_{i}=\Big\{v_{i}(\cdot)\;|
\;  v_{i}(\cdot): [0, T]\times\Omega\longrightarrow U_{i},
\hbox{ is a }\ \mathscr{G}^{i}_{t}\hbox{-adapted}\\
 \hbox{ process satisfying } \mathbb{E}\int_{0}^{T}v_{i}(t)^{2}dt<\infty\Big \},
\end{multline}
where $U_{i}$ is a nonempty convex subset of $\mathbb{R}^{r_{i}}$, and $\mathscr{G}^{i}_{t}\subseteq \mathscr{F}_{t}$ denotes
the information available to Player $i.$
Every element of $\mathscr{U}_{i}$ is called an open-loop admissible control for Player $i$ on [0, T] $ (i=1, 2)$. And
$\mathscr{U}_{1}\times\mathscr{U}_{2}$ is called the set of open-loop admissible controls for the players.  Unless otherwise stated,
Player $1$ controls $v_{1}$ and Player $2$ controls $v_{2}$.

In the following, we consider the controlled state equation of mean-field type
 {\setlength\arraycolsep{1pt}
\begin{equation}\label{Eq1}
\left\{\begin{aligned}
d&x^{v_{1},v_{2}}(t)= f\left(t, x^{v_{1},v_{2}}(t), \mathbb{E}x^{v_{1},v_{2}}(t), v_{1}(t), v_{2}(t)\right)dt\\
&+\sigma_{1}\left(t,
x^{v_{1},v_{2}}(t), \mathbb{E}x^{v_{1},v_{2}}(t), v_{1}(t), v_{2}(t)\right)dw_{1}(t)\\
&+\sigma_{2}\left(t,
x^{v_{1},v_{2}}(t), \mathbb{E}x^{v_{1},v_{2}}(t), v_{1}(t), v_{2}(t)\right)dw_{2}(t),\\
&x^{v_{1},v_{2}}(0)=x_{0},\; t\geq 0,
\end{aligned}
\right.
\end{equation}
}
and the cost functional
\begin{multline}\label{Eq4'}
J_{i}(v_{1}(\cdot), v_{2}(\cdot))\\
=\mathbb{E}\left[\int_{0}^{T}l_{i}
\left(t, x^{v_{1},v_{2}}(t), \mathbb{E}x^{v_{1},v_{2}}(t), v_{1}(t), v_{2}(t)\right)dt\right.\\
\left.+\varphi_{i}(x^{v_{1},v_{2}}(T), \mathbb{E}x^{v_{1},v_{2}}(T))\vphantom{\int_{0}^{T}}\right],
\end{multline}
where the mappings
\begin{equation*}
    \begin{aligned}
&f(t, x, \tilde{x}, v_{1}, v_{2}): \Omega\times[0,T]\times\mathbb{R}^{n+n}\times
    U_{1}\times U_{2}\rightarrow\mathbb{R}^{n},\\
&\sigma_{1}(t, x, \tilde{x}, v_{1}, v_{2}): \Omega\times[0,T]\times\mathbb{R}^{n+n}\times
    U_{1}\times U_{2}\rightarrow\mathbb{R}^{n\times d_{1}},\\
&\sigma_{2}(t, x, \tilde{x}, v_{1}, v_{2}): \Omega\times[0,T]\times\mathbb{R}^{n+n}\times
    U_{1}\times U_{2}\rightarrow\mathbb{R}^{n\times d_{2}},\\
&l_{i}(t, x, \tilde{x}, v_{1}, v_{2}): \Omega\times[0,T]\times\mathbb{R}^{n+n}\times
    U_{1}\times U_{2}\rightarrow\mathbb{R},\\
&\varphi_{i}(x, \tilde{x}): \Omega\times\mathbb{R}^{n+n}\rightarrow\mathbb{R},
    \end{aligned}
\end{equation*}
satisfy the following assumptions:
\begin{itemize}
  \item [(A1)] the coefficients $f, \sigma_{1}$ and $\sigma_{2}$ are $\mathscr{F}_{t}$-adapted and bounded by
          $C(1+|x|+|\tilde{x}|+|v_{1}|+|v_{2}|)$. They are also continuously differentiable
          with respect to $(x, \tilde{x}, v_{1}, v_{2})$ and their  partial derivatives are Lipschitz continuous and uniformly bounded.
  \item [(A2)] $l_{1}$ and $l_{2}$  are $\mathscr{F}_{t}$-adapted and continuously differentiable with respect
to $(x, \tilde{x}, v_{1}, v_{2})$. $\varphi_{1}$ and $\varphi_{2}$ are  $\mathscr{F}_{T}$-measurable and continuously differentiable with respect
to $(x, \tilde{x})$. Moreover, their partial derivatives are Lipschitz continuous and  bounded by $C(1+|x|+|\tilde{x}|+|v_{1}|+|v_{2}|)$.
\end{itemize}

Our aim is to find $(u_{1}, u_{2})\in\mathscr{U}_{1}\times\mathscr{U}_{2} $ such that
\begin{equation}\label{Eq5}
\left\{
\begin{aligned}
 J_{1}(u_{1}(\cdot), u_{2}(\cdot))\leq J_{1}(v_{1}(\cdot), u_{2}(\cdot)),\\
J_{2}(u_{1}(\cdot), u_{2}(\cdot))\leq J_{2}(u_{1}(\cdot),
v_{2}(\cdot)),
\end{aligned} \right.
\end{equation}
for all $(v_{1}, v_{2})\in\mathscr{U}_{1}\times\mathscr{U}_{2} $. We call $(u_{1}, u_{2})$ an
open-loop Nash equilibrium point of the game problem (if it exists).

Since $\mathscr{G}^{1}_{t}$ and $\mathscr{G}^{2}_{t}$ are the sub-information of $\mathscr{F}_{t}$, it implies this is the partial information game problem. On the contrary, when  $\mathscr{G}^{1}_{t}=\mathscr{G}^{2}_{t}=\mathscr{F}_{t}, t\in[0, T],$ it reduces to be a complete information case. So the problem \eqref{Ad}-\eqref{Eq5} denotes the partial information nonzero-sum differential game problem of the mean-field-type SDEs. For simplicity, we denoted it  by  Problem (MF).

\section {Nash Equilibrium Point}
\subsection{Necessary Conditions}
%......................................................................
In this subsection, we establish a necessary conditions for Nash equilibrium points of Problem (MF).
Let us suppose now that $\big(u_{1}(\cdot), u_{2}(\cdot)\big)$ is an equilibrium point with the corresponding optimal state $x(\cdot)$. Then we define the perturbed control as follows:
\begin{equation}\label{pc}
  u_{\varepsilon_{i}}(t)=u_{i}(t)+\varepsilon_{i}(v_{i}(t)-u_{i}(t)),
\end{equation}
where $\varepsilon_{i}>0$ is sufficiently small and  $v_{i}(t)$ is an arbitrary admissible control of Player $i\; (i=1,2)$.
Notice that
$U_{i}$ is convex, then for $0\leq\varepsilon_{i}\leq 1, 0\leq
t\leq T$, it yields $u_{\varepsilon_{i}}(t)\in\mathscr{U}_{i}.$  We denote by
$x_{\varepsilon_{1}}(\cdot)$ \big(resp. $x_{\varepsilon_{2}}(\cdot)$\big) the state $x^{u_{\varepsilon_{1}}, u_{2}}$ \big(resp. $x^{u_{1}, u_{\varepsilon_{2}}}$\big) associated with $\big(u_{\varepsilon_{1}}(\cdot),  u_{2}(\cdot)\big)$ \big(resp. $(u_{1}(\cdot),  u_{\varepsilon_{2}}(\cdot))$\big).
For simplicity, we set $g(t)=g\left(t, x(t), \mathbb{E}x(t), u_{1}(t), u_{2}(t)\right), \; g=f, \sigma_{1}, \sigma_{2}, l_{1}, l_{2}.$

We introduce the following variational equations:
\begin{equation}\label{Eq7}
\left\{
\begin{aligned}
 dx^{i}(t) &=\big[f_x(t)x^{i}(t)+f_{\tilde{x}}(t)\mathbb{E}x^{i}(t)\\
           &\hspace{6mm}+f_{v_{i}}(t)(v_{i}(t)-u_{i}(t))\big]dt\\
           &+\big[\sigma_{1x}(t)x^{i}(t)+\sigma_{1\tilde{x}}(t)\mathbb{E}x^{i}(t)\\
           &\hspace{6mm}+\sigma_{1v_{i}}(t)(v_{i}(t)-u_{i}(t))\big]dw_{1}(t)\\
           &+\big[\sigma_{2x}(t)x^{i}(t)+\sigma_{2\tilde{x}}(t)\mathbb{E}x^{i}(t)\\
           &\hspace{6mm}+\sigma_{2v_{i}}(t)(v_{i}(t)-u_{i}(t))\big]dw_{2}(t),\\
 x^{i}(0)& =0,\quad i=1, 2.
\end{aligned}
\right.
\end{equation}
For $i=1, 2$, we set
\begin{align*}
\bar{x}_{\varepsilon_{i}}(t)=\varepsilon&^{-1}_{i}\big(x_{\varepsilon_{i}}(t)-x(t)\big)-x^{i}(t),\\
\psi_{\varepsilon_{1}}(t)=\big(&x_{\varepsilon_{1}}(t), \mathbb{E}x_{\varepsilon_{1}}(t), u_{\varepsilon_{1}}(t), u_{2}(t)\big),\\
\psi^{\lambda}_{\varepsilon_{1}}(t)=\Big(&x(t)+\lambda\varepsilon_{1}\big(x^{1}(t)+\bar{x}_{\varepsilon_{1}}(t)\big), \mathbb{E}x(t)\\
&\hspace{20mm}+\lambda\varepsilon_{1}\mathbb{E}\big(x^{1}(t)+\bar{x}_{\varepsilon_{1}}(t)\big)\Big),\\
 \phi^{\lambda}_{\varepsilon_{1}}(t)=\Big(&x(t)+\lambda\varepsilon_{1}\big(x^{1}(t)+\bar{x}_{\varepsilon_{1}}(t)\big), \mathbb{E}x(t)\\
 &\hspace{15mm}+\lambda\varepsilon_{1}\mathbb{E}\big(x^{1}(t)+\bar{x}_{\varepsilon_{1}}(t)\big), u_{1}(t)\\
 &\hspace{20mm}+\lambda\varepsilon_{1}\big(v^{1}(t)-u_{1}(t)\big), u_{2}(t)\Big).
\end{align*}

\noindent Then by a similar method as shown in Li \cite{Li2012} and Hui and Xiao \cite{HX2012} with a minor modification, we have the following convergence result.

\begin{lemma}\label{Lemma3.2}
 Under Assumption (A1), we have
\begin{equation}\label{Eq13}
\begin{aligned}
 \lim\limits_{\varepsilon_{1}\rightarrow 0}{\mathbb{E}}\sup\limits_{0\leq t\leq T}|
 x_{\varepsilon_{1}}(t)-x(t)|^2=0,\\
 \lim\limits_{\varepsilon_{2}\rightarrow 0}{\mathbb{E}}\sup\limits_{0\leq t\leq T}|
 x_{\varepsilon_{2}}(t)-x(t)|^2=0.
\end{aligned}
\end{equation}
\end{lemma}
\emph{Proof}. By Assumption (A1) and the Burkholder-Davis-Gundy inequality, we derive
\begin{align*}
&\mathbb{E}\sup\limits_{0\leq t\leq T}|x_{\varepsilon_{1}}(t)-x(t)|^2\\
&\leq 3T\mathbb{E}\int_{0}^{T}|f\big(t, \psi_{\varepsilon_{1}}(t)\big)-f(t)|^{2}dt\\
& +12\mathbb{E}\int_{0}^{T}|\sigma_{1}\big(t, \psi_{\varepsilon_{1}}(t)\big)-\sigma_{1}(t)|^{2}dt\\
& +12\mathbb{E}\int_{0}^{T}|\sigma_{2}\big(t, \psi_{\varepsilon_{1}}(t)\big)-\sigma_{2}(t)|^{2}dt\\
&\leq C_{T}\mathbb{E}\int_{0}^{T}|x_{\varepsilon_{1}}(t)-x(t)|^2dt\\
&\hspace{4mm}+\varepsilon_{1}^{2}C_{T}\mathbb{E}\int_{0}^{T}|v_{1}(t)-u_{1}(t)|^{2}dt,
\end{align*}
where $C_{T}>0$ is a constant only depending on $T>0$ and the Lipschitz coefficients of $f$, $\sigma_{1}$ and $\sigma_{2}.$ From Gronwall's inequality we get the desired result.  \hfill$\Box$

\begin{lemma}\label{Lemma3.1} Under Assumption (A1), it yields
\begin{equation}\label{Eq9}
\begin{aligned}
 \lim\limits_{\varepsilon_{1}\rightarrow 0}{\mathbb{E}}\sup\limits_{0\leq t\leq T}|
 \bar{x}_{\varepsilon_{1}}(t)|^2=0,\\
 \lim\limits_{\varepsilon_{2}\rightarrow 0}{\mathbb{E}}\sup\limits_{0\leq t\leq T}|
 \bar{x}_{\varepsilon_{2}}(t)|^2=0.
\end{aligned}
\end{equation}
\end{lemma}
{\it Proof}. Without loss of generality, we prove the first result of \eqref{Eq9} and  the latter one can be similarly derived.
For $g=f, \sigma_{1}, \sigma_{2},$ we set
\begin{align*}
  &a^{1}_{g}(t)=\int_{0}^{1}g_{x}(t, \phi^{\lambda}_{\varepsilon_{1}}(t))d\lambda, \\
   &a^{2}_{g}(t)=\int_{0}^{1}g_{\tilde{x}}(t, \phi^{\lambda}_{\varepsilon_{1}}(t))d\lambda,\\
  &a^{3}_{g}(t)=\int_{0}^{1}\big(g_{x}(t, \phi^{\lambda}_{\varepsilon_{1}}(t))-g_{x}(t)\big)d\lambda\cdot x^{1}(t)\\
  &\hspace{5mm}+\int_{0}^{1}\big(g_{\tilde{x}}(t, \phi^{\lambda}_{\varepsilon_{1}}(t))-g_{\tilde{x}}(t)\big)d\lambda\cdot\mathbb{E}x^{1}(t)\\
  &\hspace{5mm}+\int_{0}^{1}\left(g_{v_{1}}(t, \phi^{\lambda}_{\varepsilon_{1}}(t))-g_{v_{1}}(t)\right)d\lambda\cdot\big(v_{1}(t)-u_{1}(t)\big).
\end{align*}
Due to Assumption (A1), $a^{1}_{g}$ and $a^{2}_{g}$ are both uniformly bounded and $\lim\limits_{\varepsilon_{1}\rightarrow 0}{\mathbb{E}}\left(\sup\limits_{0\leq t\leq T}|a^{3}_{g}(t)|^{2}\right)=0.$
Then we have
\begin{equation*}
\left\{
\begin{aligned}
 d\bar{x}_{\varepsilon_{1}}&(t)=\mbox{\ }[a^{1}_{f}(t)\bar{x}_{\varepsilon_{1}}(t)+a^{2}_{f}(t)\mathbb{E}\bar{x}_{\varepsilon_{1}}(t)+a^{3}_{f}(t)]dt\\
  &+[a^{1}_{\sigma_{1}}(t)\bar{x}_{\varepsilon_{1}}(t)+a^{2}_{\sigma_{1}}(t)\mathbb{E}\bar{x}_{\varepsilon_{1}}(t)+a^{3}_{\sigma_{1}}(t)]dw_{1}(t)\\
  &+[a^{1}_{\sigma_{2}}(t)\bar{x}_{\varepsilon_{1}}(t)+a^{2}_{\sigma_{2}}(t)\mathbb{E}\bar{x}_{\varepsilon_{1}}(t)+a^{3}_{\sigma_{2}}(t)]dw_{2}(t),\\
  \bar{x}_{\varepsilon_{1}}&(0)=\mbox{\ }0.
\end{aligned}
\right.
\end{equation*}

Applying It\^{o}'s formula to $|\bar{x}_{\varepsilon_{1}}(t)|^2$ and Assumption (A1), we have
\begin{equation*}
\begin{aligned}
\mathbb{E}&\left[\sup\limits_{0\leq t\leq T}|\bar{x}_{\varepsilon_{1}}(t)|^2\right]
\leq C \mathbb{E}\int_0^T|\bar{x}_{\varepsilon_{1}}(t)|^2dt\\
&+C\mathbb{E}\left[\sup\limits_{0\leq t\leq T}\left(|a_{f}^{3}(t)|^2+|a_{\sigma}^{3}(t)|^2+|a_{\sigma_{2}}^{3}(t)|^2\right)\right].
\end{aligned}
\end{equation*}
Then we can get the first convergence result of (\ref{Eq9}) from
Gronwall's inequality.          \hfill$\Box$

\par
Since $(u_{1}(\cdot), u_{2}(\cdot))$ is the Nash equilibrium point, then
it follows that
\begin{equation}\label{Eq11}
\varepsilon_{1}^{-1}[J_{1}(u_{\varepsilon_{1}}(\cdot),u_{2}(\cdot))-J_{1}(u_{1}(\cdot),
u_{2}(\cdot))]\geq 0
\end{equation}
and
\begin{equation}\label{Eq11'}
\varepsilon_{2}^{-1}[J_{2}(u_{1}(\cdot),u_{\varepsilon_{2}}(\cdot))-J_{2}(u_{1}(\cdot),
u_{2}(\cdot))]\geq 0.
\end{equation}

\begin{lemma}\label{Lemma2.3}
 Let Assumptions (A1) and (A2) hold. Then the
following {\it variational inequality} holds for i=1, 2:
\begin{multline}\label{Eq12}
\mathbb{E}\int_0^T\Big[l_{ix}(t)x^{i}(t)+l_{i\tilde{x}}(t)\mathbb{E}x^{i}(t)+l_{iv_{i}}(t)(v_{i}(t)-u_{i}(t))\Big]dt\\
  +{\mathbb{E}}[\varphi_{ix}(x(T), \mathbb{E}x(T))x^{i}(T)+\varphi_{i\tilde{x}}(x(T), \mathbb{E}x(T))\mathbb{E}x^{i}(T)]\\
  \geq 0.
\end{multline}
\end{lemma}
{\it Proof.} We firstly prove \eqref{Eq12} holds for
$i=1$ and the another case can be similarly derived. From \eqref{Eq11}, it yields
\begin{align*}
  0 \leq & \varepsilon_{1}^{-1}[J_{1}(u_{\varepsilon_{1}}(\cdot),u_{2}(\cdot))-J_{1}(u_{1}(\cdot), u_{2}(\cdot))]\\
  = & \varepsilon_{1}^{-1}\mathbb{E}\int_0^T\big[l_{1}(t, \psi_{\varepsilon_{1}}(t))-l_{1}(t)\big]dt\\
  & + \varepsilon_{1}^{-1}{\mathbb{E}}\big[\varphi_{1}(x_{\varepsilon_{1}}(T), \mathbb{E}x_{\varepsilon_{1}}(T))-\varphi_{1}(x(T), \mathbb{E}x(T))\big]\\
  = & I_{1}+I_{2}.
\end{align*}
From Assumptions (A1), (A2) and Lemma \ref{Lemma3.1}, we derive
\begin{multline}\label{Eq59}
I_{1}=\mathbb{E}\int_0^T\left[\int_0^1 l_{1x}(t, \phi^{\lambda}_{\varepsilon_{1}}(t))d\lambda\cdot\big(x^{1}(t)+\bar{x}_{\varepsilon_{1}}(t)\big)\right.\\
+\int_0^1 l_{1\tilde{x}}(t, \phi^{\lambda}_{\varepsilon_{1}}(t))d\lambda\cdot\mathbb{E}\big(x^{1}(t)+\bar{x}_{\varepsilon_{1}}(t)\big)\\
\left.+\int_0^1 l_{1v_{1}}(t, \phi^{\lambda}_{\varepsilon_{1}}(t))d\lambda\cdot\big(v_{1}(t)-u_{1}(t)\big)\vphantom{\int_0^1}\right]dt\\
\longrightarrow\mbox{\
}{\mathbb{E}}\int_0^T\Big[l_{1x}(t)x^{1}(t)+l_{1\tilde{x}}(t)\mathbb{E}x^{1}(t)\\
 +l_{1v_{1}}(t)\big(v_{1}(t)-u_{1}(t)\big)\Big]dt,
\end{multline}
\begin{multline}\label{Eq57}
I_{2}={\mathbb{E}}\left[\int_0^1
\varphi_{1x}\Big(\psi^{\lambda}_{\varepsilon_{1}}(T)\Big)d\lambda\cdot\Big(\bar{x}_{\varepsilon_{1}}(T)+x^{1}(T)\Big)\right.\\
\left.+\int_0^1
\varphi_{1\tilde{x}}\Big(\psi^{\lambda}_{\varepsilon_{1}}(T)\Big)d\lambda\cdot\mathbb{E}\Big(\bar{x}_{\varepsilon_{1}}(T)+x^{1}(T)\Big)\right]\\
  \longrightarrow  {\mathbb{E}}\Big[\varphi_{1x}(x(T), \mathbb{E}x(T))x^{1}(T)\\
  +\varphi_{1\tilde{x}}(x(T), \mathbb{E}x(T))\mathbb{E}x^{1}(T)\Big].
\end{multline}
Combining (\ref{Eq59}) with (\ref{Eq57}), the inequality (\ref{Eq12})
follows for $i=1$.  \hfill$\Box$

Next, we define the {\it Hamiltonian function}
$H_{i}:[0,T]\times{\mathbb{R}}^n\times{\mathbb{R}}^n\times U_{1}\times U_{2}\times{\mathbb{R}}^n\times\mathbb{R}^{n\times d_{1}}\times\mathbb{R}^{n\times d_{2}}\rightarrow{\mathbb{R}}$ as follows:
\begin{equation*}
\begin{aligned}
 &H_{i}(t, a, \tilde{a}, v_{1}, v_{2}, q_{i}, k_{i1}, k_{i2})\\
\triangleq &\mbox{\ }\langle q_{i},f(a, \tilde{a}, v_{1}, v_{2})\rangle+\langle
k_{i1},\sigma_{1}(a, \tilde{a}, v_{1}, v_{2})\rangle\\
&+\langle k_{i2},\sigma_{2}(a, \tilde{a}, v_{1}, v_{2})\rangle+l_{i}(a, \tilde{a}, v_{1}, v_{2}),
\end{aligned}
\end{equation*}
and denote $H_{i}(t)=H_{i}\big(t, x(t), \mathbb{E}x(t), u_{1}(t),
u_{2}(t), q_{i}(t),$ $k_{i1}(t), k_{i2}(t)\big),$ $i=1, 2.$ Let us consider the following adjoint BSDE of mean-field type
\begin{equation}\label{Eq18}
\left\{
\begin{aligned}
-dq_{i}(t)=&\mbox{\ }\big[f_{x}^{*}(t)q_{i}(t)+\mathbb{E}\big(f_{\tilde{x}}^{*}(t)q_{i}(t)\big)
                     +\sigma_{1x}^{*}(t)k_{i1}(t)\\
& +\mathbb{E}\big(\sigma_{1\tilde{x}}^{*}(t)k_{i1}(t)\big)+\sigma_{2x}^{*}(t)k_{i2}(t)\\
& +\mathbb{E}\big(\sigma_{2\tilde{x}}^{*}(t)k_{i2}(t)\big)+l_{ix}^{*}(t)+\mathbb{E}l_{i\tilde{x}}^{*}(t)\big]dt\\
&  -k_{i1}(t)dw_{1}(t)-k_{i2}(t)dw_{2}(t),\\
  q_{i}(T)=& \varphi_{ix}^{*}\big(x(T), \mathbb{E}x(T)\big)+\mathbb{E}\varphi_{i\tilde{x}}^{*}\big(x(T), \mathbb{E}x(T)\big),
\end{aligned}
\right.
\end{equation}
which coupled with \eqref{Eq1} constitutes a forward-backward SDE (FBSDE) of mean-field type.

In the sequel, we state necessary conditions of Nash equilibrium points, i.e. stochastic maximum principle as follows:
\begin{theorem}[Maximum Principle]\label{Theorem2.1}
Suppose (A1) and (A2) hold. Let $(u_{1}(\cdot), u_{2}(\cdot))$ be a Nash
equilibrium point of Problem (MF) with the corresponding
solutions $x(\cdot)$ and $
\big(q_{i}(\cdot), k_{i1}(\cdot), k_{i2}(\cdot)\big)$  of \eqref{Eq1}
and \eqref{Eq18}.
 Then it follows that
\begin{equation}\label{Eq19}
\begin{aligned}
\mathbb{E}[H_{1v_{1}}(t)|\mathscr{G}^{1}_{t}](v_{1}-u_{1}(t))\geq 0
\end{aligned}
\end{equation}
and
\begin{equation}\label{Eq19'}
\begin{aligned}
\mathbb{E}[H_{2v_{2}}(t)|\mathscr{G}^{2}_{t}](v_{2}-u_{2}(t))\geq 0,
\end{aligned}
\end{equation}
$dtd\mathbb{P}-a.e.,$ for any $v_{1}\in U_{1}$ and $v_{2}\in U_{2}.$
\end{theorem}

{\it Proof}.  We firstly prove \eqref{Eq19}. Applying It\^{o}'s formula to $\langle
x^{1}(t),q_{1}(t)\rangle$, for any $v_{1}(\cdot)\in\mathscr{U}_{1}$ we
obtain
\begin{multline}
{\mathbb{E}}\big[\varphi_{1x}\big(x(T), \mathbb{E}x(T)\big)x^{1}(T)+\varphi_{1\tilde{x}}\big(x(T), \mathbb{E}x(T)\big)\mathbb{E}x^{1}(T)\big]\\
=\mbox{\
}{\mathbb{E}}\int_0^T\Big[-l_{1x}(t)x^{1}(t)-l_{1\tilde{x}}(t)\mathbb{E}x^{1}(t)\\
-l_{1v_{1}}(t)\big(v_{1}(t)-u_{1}(t)\big)\Big]dt\\
+{\mathbb{E}}\int_0^T H_{1v_{1}}(t)\big(v_{1}(t)-u_{1}(t)\big)dt.
\end{multline}
This together with the variational inequality (\ref{Eq12}) derives
\begin{multline}
{\mathbb{E}}\int_0^T\mathbb{E}\left[ H_{1v_{1}}(t)\big|\mathscr{G}^{1}_{t}\right]\big(v_{1}(t)-u_{1}(t)\big)dt\\
={\mathbb{E}}\int_0^T H_{1v_{1}}(t)\big(v_{1}(t)-u_{1}(t)\big)dt\geq0,
\end{multline}
for all $v_{1}(\cdot)\in \mathscr{U}_{1},$ which implies (\ref{Eq19}) holds. By the
similar method above, \eqref{Eq19'} is also true. \hfill $\Box$

\subsection{Sufficient Conditions}

In what follows, we proceed to establish the sufficient conditions of Nash equilibrium points (also called verification theorem).

\begin{theorem}[Verification Theorem]\label{Theorem3.2}
Let (A1) and (A2) hold. Let  $(u_{1}(\cdot),
u_{2}(\cdot))\in\mathscr{U}_{1}\times\mathscr{U}_{2}$ with the
corresponding solutions $x$ and $ (q_{i}, k_{i1}, k_{i2})$ to
the equations \eqref{Eq1} and \eqref{Eq18}. Suppose
$H_{i}\big(t, a, \tilde{a}, v_{1}, v_{2}, q_{i}(t), k_{i1}(t), k_{i2}(t)\big)$ and $\varphi_{i}$ are convex with respect
to $(a, \tilde{a}, v_{i})$. Moreover,
\begin{multline}\label{Eq60}
\mathbb{E}\left[H_{1}\big(t,  x(t), \mathbb{E}x(t), u_{1}(t), u_{2}(t), q_{1}(t),
k_{11}(t),\right.\\
 k_{12}(t)\big)\big|\mathscr{G}^{1}_{t}\big]=\inf_{v_{1}\in\ U_{1}}\mathbb{E}\left[ H_{1}\big(t, x(t), \mathbb{E}x(t), v_{1}, u_{2}(t),\right.\\
\left. q_{1}(t), k_{11}(t), k_{12}(t)\big|\mathscr{G}^{1}_{t}\right],
\end{multline}
\begin{multline}\label{Eq60'}
\mathbb{E}\big[H_{2}\big(t,  x(t), \mathbb{E}x(t), u_{1}(t), u_{2}(t), q_{2}(t),
k_{21}(t),\\
 k_{22}(t)\big)\big|\mathscr{G}^{2}_{t}\big]
=\inf_{v_{2}\in\ U_{2}}\mathbb{E}\big[ H_{2}\big(t, x(t), \mathbb{E}x(t), u_{1}(t), v_{2},\\
 q_{2}(t), k_{21}(t), k_{22}(t)\big|\mathscr{G}^{2}_{t}\big].
\end{multline}
Then  $(u_{1}(\cdot), u_{2}(\cdot))$ is a Nash equilibrium point of Problem (MF).
\end{theorem}
{\it Proof}. Let $v_{i}(\cdot)\in\mathscr{U}_{i}, i=1,2.$ We denote by $x^{v_{1}}$ and $x$ the solutions to \eqref{Eq1} associated with the admissible controls $(v_{1}, u_{2})$ and $(u_{1}, u_{2})$, respectively. We set
\begin{multline*}
  H_{1}(t)=H_{1}\big(t, x(t), \mathbb{E}x(t), u_{1}(t), u_{2}(t), q_{1}(t), k_{11}(t),\\
   k_{12}(t)\big), \quad H_{1}^{v_{1}}(t)=H_{1}\big(t, x^{v_{1}}(t), \mathbb{E}x^{v_{1}}(t), v_{1}(t), u_{2}(t),\\
    q_{1}(t), k_{11}(t), k_{12}(t)\big), \; g(t)=g\big(t, x(t), \mathbb{E}x(t), u_{1}(t), u_{2}(t)\big),\\
     g^{v_{1}}(t)=g\big(t, x^{v_{1}}, \mathbb{E}x^{v_{1}}(t), v_{1}(t), u_{2}(t)\big),\quad g=f, \sigma_{1}, \sigma_{2}, l_{1}.
\end{multline*}
By virtue of the convexity of $\varphi_{1} $, we have for any $v_{1}(\cdot)\in \mathscr{U}_{1}$
\begin{equation}\label{Sq1}
    \begin{aligned}
    J_{1}(v_{1}(\cdot), u_{2}(\cdot))-J_{1}(u_{1}(\cdot), u_{2}(\cdot))\geq I_{1}+I_{2}
    \end{aligned}
\end{equation}
with
\begin{align}
I_{1}&=\mathbb{E}\big[\varphi_{1x}(x(T), \mathbb{E}x(T))(x^{v_{1}}(T)-x(T))\nonumber\\
&\hspace{8mm}+\varphi_{1\tilde{x}}(x(T), \mathbb{E}x(T))\mathbb{E}(x^{v_{1}}(T)-x(T))\big],\nonumber\\
I_{2}&=\mathbb{E}\int_{0}^{T}\Big(l^{v_{1}}_{1}(t)-l_{1}(t)\Big)dt\nonumber\\
  & =\mathbb{E}\int_{0}^{T}\Big[H_{1}^{v_{1}}(t)-H_{1}(t)-\langle q_{1}(t), f^{v_{1}}(t)-f(t)\rangle\nonumber\\
       &-\langle k_{11}(t), \sigma_{1}^{v_{1}}(t)-\sigma_{1}(t)\rangle\nonumber\\
       &-\langle k_{12}(t), \sigma_{2}^{v_{1}}(t)-\sigma_{2}(t)\rangle\Big]dt.\label{Eq23}
\end{align}
Applying It\^{o}'s formula to $ \langle q_{1}(t), x^{v_{1}}(t)-x(t)\rangle,$ we have
\begin{align}
I_{1}=\ & \mathbb{E}\int_{0}^{T}\Big[\langle q_{1}(t), f^{v_{1}}(t)-f(t)\rangle\nonumber\\
&+\langle k_{11}(t), \sigma_{1}^{v_{1}}(t)-\sigma_{1}(t)\rangle+\langle k_{12}(t), \sigma_{2}^{v_{1}}(t)-\sigma_{2}(t)\rangle\nonumber\\
       &-\langle H_{1x}^{*}(t), x^{v_{1}}(t)-x(t)\rangle\nonumber\\
 \label{Eq23'}      &-\langle \mathbb{E}H_{1\tilde{x}}^{*}(t), x^{v_{1}}(t)-x(t)\rangle\Big]dt.
\end{align}
Substituting \eqref{Eq23} and \eqref{Eq23'} into \eqref{Sq1} and applying the convexity of $H_{1}$, we get
\begin{align}
   & J_{1}(v_{1}(\cdot), u_{2}(\cdot))-J_{1}(u_{1}(\cdot), u_{2}(\cdot))\nonumber\\
 \geq &\ \mathbb{E}\int_{0}^{T}\Big(H_{1}^{v_{1}}(t)-H_{1}(t)-\langle
H_{1x}^{*}(t), x^{v_{1}}(t)-x(t)\rangle\nonumber\\
&\hspace{10mm}-\langle
\mathbb{E}H_{1\tilde{x}}^{*}(t), x^{v_{1}}(t)-x(t)\rangle\Big) dt\nonumber\\
 \geq &\ \mathbb{E}\int_{0}^{T}H_{1v_{1}}(t)\big(v_{1}(t)-u_{1}(t)\big)dt\nonumber\\
 = &\  \mathbb{E}\int_{0}^{T}\mathbb{E}\left[H_{1v_{1}}(t)\big|\mathscr{G}^{1}_{t}\right]\big(v_{1}(t)-u_{1}(t)\big)dt.\label{Eq62}
   \end{align}
Condition \eqref{Eq60} implies $\mathbb{E}\left[H_{1v_{1}}(t)\big|\mathscr{G}^{1}_{t}\right]\big(v_{1}(t)-u_{1}(t)\big)\geq 0,$ $dtd\mathbb{P}-a.e.$
on $[0, T],$ which derives $$J_{1}(v_{1}(\cdot), u_{2}(\cdot))-J_{1}(u_{1}(\cdot), u_{2}(\cdot))\geq 0.$$
 Similarly, we can also derive
$J_{2}(u_{1}(\cdot), v_{2}(\cdot))-J_{2}(u_{1}(\cdot), u_{2}(\cdot))\geq 0.$ The proof is completed.        \hfill$\Box$

\section{LQ Example}
In this section, we work out an LQ example to illustrate the theoretical result. Without loss of generality, we only consider the following case: $n=d_{1}=d_{2}=1, b_{1}(t)b_{2}(t)\neq 0, t\in[0, T].$ Throughout this section, we assume additional condition.

\begin{itemize}
  \item [(A3)] $ m_{1}^{-1}(t)b_{1}^{2}(t)=m_{2}^{-1}(t)b_{2}^{2}(t), t\in[0, T].$
\end{itemize}

\begin{example}\label{Exm2}
Consider the system of linear mean-field SDE
\begin{equation}\label{Eq38}
\left\{
\begin{aligned}
dx^{v_{1}, v_{2}}(t)=& \big[a(t)x^{v_{1}, v_{2}}(t)+\bar{a}(t)\mathbb{E}x^{v_{1}, v_{2}}(t)\\
 &+b_{1}(t)v_{1}(t)+b_{2}(t)v_{2}(t)\big]dt\\
         &+c_{1}(t)dw_{1}(t)+c_{2}(t)dw_{2}(t),\\
x^{v_{1}, v_{2}}(0)=&\ x_{0},
\end{aligned}
\right.
\end{equation}
with the quadratic cost functional
\begin{multline}\label{Eq39}
J_{i}\big(v_{1}(\cdot), v_{2}(\cdot)\big)=\frac{1}{2}
\mathbb{E}\left[\int_{0}^{T}\Big(g_{i}(t)\big(x^{v_{1}, v_{2}}(t)\big)^{2}\right.\\
+\bar{g}_{i}(t)\big(\mathbb{E}x^{v_{1}, v_{2}}(t)\big)^{2}
  +m_{i}(t)\big(v_{i}(t)\big)^{2}\Big)dt\\
  \left.+h_{i}\big(x^{v_{1}, v_{2}}(T)\big)^{2}+\bar{h}_{i}\big(\mathbb{E}x^{v_{1}, v_{2}}(T)\big)^{2}\vphantom{\int_{0}^{T}}\right]
\end{multline}
and the information available to two players $\mathscr{G}^{1}_{t}=\mathscr{G}^{1}_{t}=\mathscr{F}_{t}^{w_{1}}=\sigma\{w_{1}(s), 0\leq s\leq t\}.$

Here, all coefficients with respect to $t$ in \eqref{Eq38} and \eqref{Eq39} are deterministic and uniformly bounded. In addition, $g_{i}$ and $\bar{g}_{i}$ are non-negative, $h_{i}$ and $\bar{h}_{i}$ are non-negative constants, and $m_{i}$ is positive. The set of admissible controls for Player $i$ is defined by
\begin{multline}\label{Eq40}
    \mathscr{U}_{i}=\{v_{i}(\cdot)\;|\;  v_{i}(\cdot) \hbox{ is an } \mathbb{R}\hbox{-valued }\mathscr{F}_{t}^{w_{1}}\hbox{-adapted process }\\ \hbox{ satisfying }\mathbb{E}\int_{0}^{T}v_{i}^{2}(t)dt<\infty \}, i=1, 2.
 \end{multline}
Then the unique Nash equilibrium point is denoted by
\begin{equation}\label{Eq21}
  \left\{
  \begin{aligned}
    & u_{1}(t)=-m_{1}^{-1}(t)b_{1}(t)\big(\tau_{1}(t)\hat{x}(t)+\delta_{1}(t)\mathbb{E}x(t)\big),\\
    & u_{2}(t)=-m_{2}^{-1}(t)b_{2}(t)\big(\tau_{2}(t)\hat{x}(t)+\delta_{2}(t)\mathbb{E}x(t)\big),
  \end{aligned}
 \right.
\end{equation}
where $\mathbb{E}x$, $(\tau_{1}, \tau_{2}),$ $(\delta_{1}, \delta_{2})$ and $\hat{x}$ are determined by \eqref{Eq14}, \eqref{Eq16}, \eqref{Eq17} and \eqref{Eq20}, respectively.
\end{example}

\emph{Proof}. We shall complete the proof by two parts.

\emph{\textbf{Part 1.}} We first need to prove the unique Nash equilibrium point can be represented by
\begin{equation}\label{Eq2}
  \left\{
  \begin{aligned}
    & u_{1}(t)=-m_{1}^{-1}(t)b_{1}(t)\hat{q}_{1}(t),\\
    & u_{2}(t)=-m_{2}^{-1}(t)b_{2}(t)\hat{q}_{2}(t),
  \end{aligned}
 \right.
\end{equation}
where $x$, $\big(q_{1}, k_{11}, k_{12}\big)$ and $\big(q_{2}, k_{21}, k_{22}\big)$ are the unique solution of the following coupled FBSDE of mean-field type
\begin{equation}\label{Eq3a}
\left\{
\begin{aligned}
d&x(t)= \big[a(t)x(t)+\bar{a}(t)\mathbb{E}x(t)-m_{1}^{-1}(t)b_{1}^{2}(t)\hat{q}_{1}(t)\\
&-m_{2}^{-1}(t)b_{2}^{2}(t)\hat{q}_{2}(t)\big]dt+c_{1}(t)dw_{1}(t)+c_{2}(t)dw_{2}(t),\\
x&(0)=x_{0},\\
\end{aligned}\right.
\end{equation}
\begin{equation}\label{Eq3b}
\left\{
\begin{aligned}
-&dq_{1}(t)=\big[a(t)q_{1}(t)+\bar{a}(t)\mathbb{E}q_{1}(t)+g_{1}(t)x(t)\\
           &+\bar{g}_{1}(t)\mathds{E}x(t)\big]dt-k_{11}(t)dw_{1}(t)-k_{12}(t)dw_{2}(t),\\
q&_{1}(T)=h_{1}x(T)+\bar{h}_{1}\mathbb{E}x(T),\\
\end{aligned}\right.
\end{equation}
and
\begin{equation}\label{Eq3c}
\left\{\begin{aligned}
-&dq_{2}(t)=\big[a(t)q_{2}(t)+\bar{a}(t)\mathbb{E}q_{2}(t)+g_{2}(t)x(t)\\
           &+\bar{g}_{2}(t)\mathds{E}x(t)\big]dt-k_{21}(t)dw_{1}(t)-k_{22}(t)dw_{2}(t),\\
q&_{2}(T)=h_{2}x(T)+\bar{h}_{2}\mathbb{E}x(T).
\end{aligned}\right.
\end{equation}
Here we denote by $\hat{p}(t)$ the mathematical expectation of $p(t)$ with respect to $\mathscr{F}^{w_{1}}_{t},$ i.e., $\hat{p}(t)\triangleq\mathbb{E}\left[p(t)\big|\mathscr{F}^{w_{1}}_{t}\right], p=q_{1}, q_{2}, k_{11}, k_{21}, x.$
The rest of Part 1 is divided into the following two steps.

\emph{\textbf{Step (i)}} $(u_{1}, u_{2})$ of the form \eqref{Eq2} is the Nash equilibrium point indeed.

We first write down the Hamiltonian function
\begin{align}
 &H_{i}\big(t, x, \tilde{x}, v_{1}, v_{2},  q_{i}, k_{i1}, k_{i2}\big)\nonumber\\
\triangleq&\mbox{\ } q_{i}\big(a(t)x+\bar{a}(t)\tilde{x}+b_{1}(t)v_{1}+b_{2}(t)v_{2}\big)+ k_{i1}c_{1}\nonumber\\
&+k_{i2}(t)c_{2}+\frac{1}{2}\big(g_{i}x^{2}+\bar{g}_{i}\tilde{x}^{2}+m_{i}v_{i}^{2}\big).\label{Eq4}
\end{align}
Applying Theorem \ref{Theorem2.1}, we derive the candidate Nash equilibrium point of the form \eqref{Eq2} and the coupled FBSDE of mean-field type \eqref{Eq3a}-\eqref{Eq3c}.
We can check that  $\varphi_{i}(x, \tilde{x})=\frac{1}{2}(h_{i}x^{2}+\bar{h}_{i}\tilde{x}^{2})$ and $H_{i}\big(t, x, \tilde{x}, v_{1}, v_{2},  q_{i}, k_{i1}, k_{i2}\big)$ in \eqref{Eq4} satisfy the conditions in Theorem \ref{Theorem3.2}. Therefore, $\big(u_{1}(\cdot), u_{2}(\cdot)\big)$ of the form \eqref{Eq2} is the Nash equilibrium point indeed.

Based on the arguments in Step (i), we conclude that the existence and uniqueness of the Nash equilibrium points are equivalent to the existence and uniqueness of the solutions to \eqref{Eq3a}-\eqref{Eq3c}.

\emph{\textbf{Step (ii)}} The solutions of \eqref{Eq3a}-\eqref{Eq3c} are existent and unique.

Taking mathematical expectation on both sides of \eqref{Eq3a}-\eqref{Eq3c}, we have the following forward-backward ordinary equations
\begin{equation}\label{Eq4a}
\left\{
\begin{aligned}
&d\mathbb{E}x(t)= \big[(a(t)+\bar{a}(t))\mathbb{E}x(t)\\
&\;\; -m_{1}^{-1}(t)b_{1}^{2}(t)\mathbb{E}q_{1}(t)-m_{2}^{-1}(t)b_{2}^{2}(t)\mathbb{E}q_{2}(t)\big]dt,\\
&\mathbb{E}x(0)=x_{0},\\
\end{aligned}\right.
\end{equation}
\begin{equation}\label{Eq4b}
\left\{
\begin{aligned}
-d\mathbb{E}q_{1}(t)=&\big[(a(t)+\bar{a}(t))\mathbb{E}q_{1}(t)\\
           &+(g_{1}(t)+\bar{g}_{1}(t))\mathds{E}x(t)\big]dt,\\
\mathbb{E}q_{1}(T)=&(h_{1}+\bar{h}_{1})\mathbb{E}x(T),\\
\end{aligned}\right.
\end{equation}
and
\begin{equation}\label{Eq4c}
\left\{\begin{aligned}
-d\mathbb{E}q_{2}(t)=&\big[(a(t)+\bar{a}(t))\mathbb{E}q_{2}(t)\\
           &+(g_{2}(t)+\bar{g}_{2}(t))\mathds{E}x(t)\big]dt,\\
\mathbb{E}q_{2}(T)=&(h_{2}+\bar{h}_{2})\mathbb{E}x(T).
\end{aligned}\right.
\end{equation}
Applying the method as shown in Chang and Xiao \cite{CX2014} and Assumption (A3), we can prove there exists a unique solution $(\mathbb{E}x(t), \mathbb{E}q_{1}(t), \mathbb{E}q_{2}(t))$ to \eqref{Eq4a}-\eqref{Eq4c} with the relations as follows:
\begin{equation}\label{Eq8}
  \mathbb{E}q_{i}(t)=\alpha_{i}(t)\mathbb{E}x(t),\;\; i=1, 2,
\end{equation}
\begin{multline}\label{Eq14}
  \mathbb{E}x(t)=x_{0}e^{\int_{0}^{t}[a(s)+\bar{a}(s)-m_{1}^{-1}(s)b_{1}^{2}(s)(\alpha_{1}(s)+\alpha_{2}(s))]ds},
\end{multline}
where $(\alpha_{1}, \alpha_{2})$ is the unique solution of the following Riccati equations
\begin{multline}\label{Eq6a}
\dot{\alpha_{1}}+2(a+\bar{a})\alpha_{1}-m_{1}^{-1}b_{1}^{2}\alpha_{1}^{2}-m_{2}^{-1}b_{2}^{2}\alpha_{1}\alpha_{2}\\
+g_{1}+\bar{g}_{1}=0,
\end{multline}
\begin{multline}\label{Eq6b}
\dot{\alpha_{2}}+2(a+\bar{a})\alpha_{2}-m_{1}^{-1}b_{1}^{2}\alpha_{1}\alpha_{2}-m_{2}^{-1}b_{2}^{2}\alpha_{2}^{2}\\
+g_{2}+\bar{g}_{2}=0,
\end{multline}
subject to $\alpha_{i}(T)=h_{i}+\bar{h}_{i}.$

Substituting \eqref{Eq8} into \eqref{Eq3b} and \eqref{Eq3c}, taking conditional mathematical expectation on both sides of \eqref{Eq3a}-\eqref{Eq3c} with respect to $\mathscr{F}^{w_{1}}_{t}$ and applying Lemma 5.4 in Xiong \cite{Xiong}, we have
\begin{equation}\label{Eq10a}
\left\{
\begin{aligned}
d&\hat{x}(t)= \big[a(t)\hat{x}(t)+\bar{a}(t)\mathbb{E}x(t)-m_{1}^{-1}(t)b_{1}^{2}(t)\hat{q}_{1}(t)\\
&-m_{2}^{-1}(t)b_{2}^{2}(t)\hat{q}_{2}(t)\big]dt+c_{1}(t)dw_{1}(t),\\
\hat{x}&(0)=x_{0},\\
\end{aligned}\right.
\end{equation}
\begin{equation}\label{Eq10b}
\left\{
\begin{aligned}
&-d\hat{q}_{1}(t)=\big[a(t)\hat{q}_{1}(t)+g_{1}(t)\hat{x}(t)\\
           &\; +\big(\bar{a}(t)\alpha_{1}(t)+\bar{g}_{1}(t)\big)\mathds{E}x(t)\big]dt-\hat{k}_{11}(t)dw_{1}(t),\\
&\hat{q}_{1}(T)=h_{1}\hat{x}(T)+\bar{h}_{1}\mathbb{E}x(T),\\
\end{aligned}\right.
\end{equation}
\begin{equation}\label{Eq10c}
\left\{\begin{aligned}
&-d\hat{q}_{2}(t)=\big[a(t)\hat{q}_{2}(t)+g_{2}(t)\hat{x}(t)\\
           &\; +\big(\bar{a}(t)\alpha_{2}(t)+\bar{g}_{2}(t)\big)\mathds{E}x(t)\big]dt-\hat{k}_{21}(t)dw_{1}(t),\\
&\hat{q}_{2}(T)=h_{2}\hat{x}(T)+\bar{h}_{2}\mathbb{E}x(T),
\end{aligned}\right.
\end{equation}
which constitute a kind of fully coupled forward-backward stochastic filtering equations of mean-field type and exist the unique solution $(\hat{x}, \hat{q}_{1}, \hat{k}_{11}, \hat{q}_{2}, \hat{k}_{21})$.

Similarly, based on the known $\mathbb{E}x, \mathbb{E}q_{1}, \mathbb{E}q_{2}, \hat{q}_{1}$ and $\hat{q}_{2}$, there also exists the unique
solution $\big(x, (q_{1}, k_{11}, k_{12}), (q_{2}, k_{21}, k_{22})\big)$ to \eqref{Eq3a}-\eqref{Eq3c}.

\emph{\textbf{Part 2.}} We need to verify the feedback form of the Nash equilibrium point in \eqref{Eq2} is represented by \eqref{Eq21}.

Based on the terminal conditions in  \eqref{Eq10b} and \eqref{Eq10c}, we set
\begin{equation}\label{Eq15}
  \hat{q}_{i}(t)=\tau_{i}(t)\hat{x}(t)+\delta_{i}(t)\mathbb{E}x(t), \quad i=1, 2,
\end{equation}
subject to $\tau_{i}(T)=h_{i}$ and $\delta_{i}(T)=\bar{h}_{i}$.
Applying It\^{o}'s formula to $\hat{q}_{1}(t)$ \big(resp. $\hat{q}_{2}(t)$\big) in \eqref{Eq15} and comparing the coefficients of $\hat{x}(t)$ and $\mathbb{E}x(t)$ between it and \eqref{Eq10b} \big(resp. \eqref{Eq10c}\big), respectively, we get
\begin{equation}\label{Eq16}
  \left\{
  \begin{aligned}
    \dot{\tau_{1}}+2a\tau_{1}-m_{1}^{-1}b_{1}^{2}\tau_{1}^{2}-m_{2}^{-1}b_{2}^{2}\tau_{1}\tau_{2}+g_{1}=0, \\
    \dot{\tau_{2}}+2a\tau_{2}-m_{2}^{-1}b_{2}^{2}\tau_{2}^{2}-m_{1}^{-1}b_{1}^{2}\tau_{1}\tau_{2}+g_{2}=0,
  \end{aligned}
  \right.
\end{equation}
\begin{equation}\label{Eq17}
  \left\{
  \begin{aligned}
   & \dot{\delta_{1}}+\big(2a+\bar{a}-m_{1}^{-1}b_{1}^{2}\alpha_{1}-m_{2}^{-1}b_{2}^{2}\alpha_{2}-m_{1}^{-1}b_{1}^{2}\tau_{1}\big)\delta_{1}\\   &\hspace{25mm}-m_{2}^{-1}b_{2}^{2}\tau_{1}\delta_{2}+\bar{a}\tau_{1}+\bar{a}\alpha_{1}+\bar{g}_{1}=0, \\
    & \dot{\delta_{2}}+\big(2a+\bar{a}-m_{1}^{-1}b_{1}^{2}\alpha_{1}-m_{2}^{-1}b_{2}^{2}\alpha_{2}-m_{2}^{-1}b_{2}^{2}\tau_{2}\big)\delta_{2}\\   &\hspace{25mm}-m_{1}^{-1}b_{1}^{2}\tau_{2}\delta_{1}+ \bar{a}\tau_{2}+\bar{a}\alpha_{2}+\bar{g}_{2}=0, \\
  \end{aligned}
  \right.
\end{equation}
with $\tau_{i}(T)=h_{i}$ and $\delta_{i}(T)=\bar{h}_{i}.$
Applying the method as shown in Chang and Xiao \cite{CX2014} and Assumption (A3), there exist the unique solutions to \eqref{Eq16} and \eqref{Eq17}.

Substituting \eqref{Eq15} into \eqref{Eq10a}, we can derive the explicit solution as follow:
\begin{align}
& \hat{x}(t)=x_{0}\Phi_{0}^{t} +\int_{0}^{t}\Phi_{s}^{t}c_{1}(s)dw_{1}(s)\nonumber\\
\label{Eq20} &  +\int_{0}^{t}\Phi_{s}^{t}\big(\bar{a}-m_{1}^{-1}b_{1}^{2}\delta_{1}-m_{2}^{-1}b_{2}^{2}\delta_{2}\big)(s)\mathbb{E}x(s)ds
\end{align}
with $\Phi_{s}^{t}=\exp\{\int_{s}^{t}(a-m_{1}^{-1}b_{1}^{2}\tau_{1}-m_{2}^{-1}b_{2}^{2}\tau_{2})(r)dr\}.$  The proof is completed.
\hfill $\Box$

\section{Conclusion Remarks}

In this paper, we study non-zero sum mean-field game with partial information and derive the stochastic maximum principle and verification theorem for the Nash equilibrium points.  Compared with the existing literature, the contributions of this paper are:
\begin{itemize}
  \item Partial information is more general case than complete information. The results partly generalize the related mean-field control or game problems with complete information (see e.g. \cite{BSY2013,BDL2011,Li2012,Yong2013});
  \item Mean-field-type forward-backward stochastic filtering equations are found, which enriches the theory of classical filtering;
  \item The unique Nash equilibrium point in the LQ example is represented in the feedback form of not only the optimal filtering but also the expected value of the system state, through the solutions of some Riccati equations.
\end{itemize}

In addition, since there are many partial information mean-field game problems in finance and economics, we hope the results have applications in these related areas.

\end{document}